\documentclass[11pt,reqno]{amsart}
\usepackage{amsmath,amssymb,amsthm,enumerate}
\usepackage{bm, soul, color}
\usepackage{amsfonts, mathrsfs}
\usepackage{verbatim}

\usepackage[pdftex]{graphicx}

\newcommand{\R}{\mathbb{R}}
\newcommand{\Q}{\mathbb{Q}}
\newcommand{\Z}{\mathbb{Z}}
\newcommand{\N}{\mathbb{N}}
\newcommand{\J}{\mathcal{J}}

\newtheorem{thm}{Theorem}
\newtheorem{lem}{Lemma}

\newtheorem{cor}[thm]{Corollary}
\theoremstyle{definition}

\theoremstyle{remark}

\author{Jing-Jing Huang and Jason J. Liu}
\address{JH: Department of Mathematics and Statistics, University of Nevada, Reno,
1664 N. Virginia St., Reno, NV 89557}
\email{jingjingh@unr.edu}
\address{JL: Davidson Academy, 1164 N Virginia St, Reno, NV 89503
}

\email{jliu@davidsonacademy.unr.edu}
\thanks{The first-named author is supported by the UNR VPRI startup grant 1201-121-2479}


\begin{document}

\title{Simultaneous approximation on affine subspaces}
\begin{abstract}
We solve the convergence case of the generalized Baker-Schmidt problem for simultaneous approximation on affine subspaces, under natural diophantine type conditions. In one of our theorems, we do not require monotonicity on the approximation function. In order to prove these results, we establish asymptotic formulae for the number of rational points close to an affine subspace. One key ingredient is a sharp upper bound on a certain sum of reciprocals of fractional parts associated with the matrix defining the affine subspace.
 \end{abstract}
\maketitle

\section{Diophantine approximation on affine subspaces}\label{s1}

Let $\psi:\N\rightarrow\R_{\ge0}$ be a function, which from now on we will call an \emph{approximation function}. Let $\| y\|$ be the distance of $y$ to the nearest integer for $y\in\mathbb{R}$ and $\displaystyle\|\bm{y}\|=\max_{1\le k\le n}\|y_k\|$ for a vector $\bm{y}=(y_1, y_2, \ldots, y_n)\in\R^n$. We denote
\begin{equation*}
\mathscr{S}_n(\psi):=\{\bm{x}\in\R^n:\exists^\infty{q}\in\N\textrm{ such that } \|{q}\bm{x}\|<\psi({q})\},
\end{equation*}
where $\exists^\infty$ means ``there exist infinitely many''.

 The points in the set $\mathscr{S}_n(\psi)$ are said to be simultaneously \emph{$\psi$-approximable}. A classical theorem of Khintchine \cite{Kh} states that the set $\mathscr{S}_n(\psi)$ obeys a beautiful zero-one law according to whether $\sum \psi(q)^n$ converges, which in the divergence case requires  the additional assumption that $\psi$ is non-increasing. Modern developments reveal that the problem becomes much deeper when the underlying space is restricted to a proper submanifold of $\R^n$. More precisely, we say a submanifold $\mathcal{M}$ of $\R^n$ is of \emph{Khintchine type} (resp. \emph{strong Khintchine type}) \emph{for convergence} if $|\mathcal{M}\cap\mathscr{S}_n(\psi)|_\mathcal{M}=0$ when $\sum_{n=1}^\infty\psi(q)^n$ converges and $\psi$ is a non-increasing approximation function (resp. any approximation function). Here $|\cdot|_\mathcal{M}$ means the induced Lebesgue measure on $\mathcal{M}$.

A proper affine subspace of $\R^n$ defined over $\Q$ is well known to not be of Khintchine type for convergence, due to the fact that there are too
many rational points on it. Therefore some conditions have to be imposed to obtain Khintchine type manifolds for convergence.  One possible such
assumption is that the manifold is non-degenerate, namely it deviates from any hyperplane. Indeed, some curved manifolds have been shown to be of
strong Khintchine type for convergence \cite{BVVZ, hua4, Sim}. It transpires, however, that non-degeneracy is not necessary,
something which previously had not been shown in this context. The main purpose of this paper is to address this issue for  affine subspaces, which are the most notable examples of degenerate manifolds, under some natural diophantine type conditions. 

We need to introduce some notation first. Throughout the paper, $|\bm{y}|$ denotes the supremum norm of the vector $\bm{y}$ and all vectors are regarded as row vectors by default.  Let
$$
\bm{A}\in \R^{d\times(n-d)}, \bm{\widetilde{A}}={\bm{\alpha}_0\choose\bm{A}}\in \R^{(d+1)\times(n-d)}
$$
and
 \begin{equation}\label{e1.1}
        \mathscr{H} :=\left\{(\bm{x},\widetilde{\bm{x}}\cdotp\bm{\widetilde{A}})\bigm|\bm{x}\in{\R^d} \right\},
\end{equation}
where
$$
\widetilde{\bm{x}}=(1,\bm{x}).
$$ 

For a matrix $\bm{M} \in \R^{l\times m}$, let $\mathrm{row}_u(\bm{M})$ and $\mathrm{col}_v(\bm{M})$ represent the $u$-th row and $v$-th column of $\bm{M}$. We say $\bm{M} $  is of (multiplicative) diophantine type $\omega(\bm{M})$ if $\omega(\bm{M})$ is the least upper bound of
all $\omega$ for which $$\prod_{u = 1}^{l} \|\bm{j} \cdot \mathrm{row}_u(\bm{M})\| < |\bm{j}|^{-\omega}$$ has infinitely many solutions
$\bm{j} \in \Z^{m}$. It is well known that $\omega(\bm{M})\in[m,\infty]$ and $\omega(\bm{M})=m$ for almost all matrices $\bm{M}$, with
respect to the Lebesgue measure on $\R^{l\times m}$. We refer the reader to a paper by Bugeaud \cite{Bu} for a detailed survey on 
multiplicative Diophantine approximation. In fact, using a theorem of Hussain and Simmons \cite{HS}, it is easy to see that for $\omega\ge m$ we have
$$
\dim\{M\in\R^{l\times m}:\omega(M)\ge\omega\}\le lm-1+\frac{m+1}{\omega+1}.
$$

\begin{thm}\label{t2} Let $\mathscr{H}$ be given by \eqref{e1.1}. Then
\begin{enumerate}[(a)]
\item
$\mathscr{H}$ is of strong Khintchine type for convergence if $\omega(\bm{A}) < dn$;
\item
$\mathscr{H}$ is of Khintchine type for convergence if $\omega(\widetilde{\bm{A}}) < (d+1)n$.
\end{enumerate}

\end{thm}

Previously, the only known affine subspaces of Khintchine type for convergence are lines through the origin satisfying $\omega(\bm{A}) < n$ \cite{Ko} and certain coordinate hyperplanes \cite[Theorem 5(b)]{R}. To the best of our knowledge, no affine subspace has been shown to be of strong Khintchine type for convergence until now. 

The connection between diophantine exponents of a matrix and metric diophantine properties of the affine subspace defined by the matrix is of course not new. Actually, Kleinbock  has shown that an affine subspace is extremal, i.e. $|\mathscr{H}\cap\mathscr{S}_n(q^{-\nu})|_\mathscr{H}=0$ whenever $\nu>1/n$, if and only if certain diophantine exponents $\omega_i(\widetilde{\bm{A}})$, $1\le i\le n-d$ are all less than or equal to $n$ \cite{Kl1, Kl2}. In view of the fact that the extremality theory is a special case of the Khintchine-Groshev theory, Ghosh has proved that an affine subspace $\mathscr{H}$ is of Groshev type for convergence if $\omega_i(\widetilde{\bm{A}})<n$ for all $1\le i\le n-d$ \cite{Gh}. Due to the technicality, we refer the readers to the above cited papers for the precise definition of those higher order exponents of a matrix. Here we only observe that the first exponent $\omega_1(\widetilde{\bm{A}})$ there is closely related to our $\omega(\widetilde{\bm{A}})$ in that $\omega(\widetilde{\bm{A}})\ge (d+1)\omega_1(\widetilde{\bm{A}})$. So the assumption $\omega(\widetilde{\bm{A}})<(d+1)n$ in Theorem \ref{t2} (b) is stronger than $\omega_1(\widetilde{\bm{A}})<n$. However, we need to emphasize that for the case of simultaneous approximation (Khintchine theory) considered in the current paper, we do not need any assumptions on the higher order exponents as has been done in the case of dual approximation (Groshev theory). Actually, Kleinbock has conjectured that the assumptions on the higher order exponents might be redundant in order for an affine subspace to be extremal, that is, the sole condition $\omega_1(\widetilde{\bm{A}})\le n$ suffices for $\mathscr{H}$ to be extremal \cite{Kl2}. Therefore, we may regard Theorem \ref{t2} (b) as achieving a slightly weaker form of this conjecture.

We remark in passing that there is also a divergence counterpart of the Khintchine theory for affine subspaces. This aspect has been studied for  affine coordinate subspaces, formed by fixing one or more coordinates \cite{RSS}.

In addition to the Khintchine theory, one may as well investigate the deeper Hausdorff measure of the underlying set $\mathcal{M}\cap\mathscr{S}_n(\psi)$. 
The central conjecture in this respect, is widely known as \emph{the generalized Baker-Schmidt problem}. There has been some tremendous progress towards the full solution of this problem for non-degenerate submanifolds. The divergence case is now reasonably well understood \cite{bere, BDV, BVVZ1}, and the convergence case is known for planar curves \cite{VV} and hypersurfaces \cite{hua3}. We obtain analogous results in the convergence case for affine subspaces, of which Theorem \ref{t2} is a special case.

\begin{thm}\label{t1} Let $\mathscr{H}$ be given by \eqref{e1.1} and $\psi$ be any approximation function. We have
$$\mathcal{H}^s(\mathscr{S}_n(\psi) \cap \mathscr{H}) = 0$$ when $\displaystyle\sum_{q = 1}^\infty\psi(q)^{n - d + s}q^{d - s}$ converges and
$\omega(\bm{A}) < \frac{d(n - d + s)}{d + 1 - s}.$
\end{thm}

We emphasize that in Theorem \ref{t1} the approximation function $\psi$ is not required to be monotonic. Hence, we have shown that, subject to a certain diophantine type condition, affine subspaces are of strong Jarn\'{i}k type for convergence.

\begin{thm}\label{t3} Let $\mathscr{H}$ be given by \eqref{e1.1} and $\psi$ be a non-increasing approximation function. We have
$$\mathcal{H}^s(\mathscr{S}_n(\psi) \cap \mathscr{H}) = 0$$ when $\displaystyle\sum_{q = 1}^\infty\psi(q)^{n - d + s}q^{d - s}$ converges and
$\omega(\bm{\widetilde{A}}) < \frac{(d+1)(n - d + s)}{d + 1 - s}$.
\end{thm}

The diophantine type condition in Theorem \ref{t3} is less restrictive than that in Theorem \ref{t1}, since we can  take advantage of the diophantine type information of   $\bm{\widetilde{A}}$ in lieu of $\bm{A}$ by averaging over $q$.

As the typical values of both $\omega(\bm{A})$ and $\omega(\bm{\widetilde{A}})$ are  $n-d$, Theorems \ref{t1} and \ref{t3} hold for generic affine subspace  $\mathscr{H}$ when $s>1-\frac{d}n$ and $s>0$ respectively. Compared to the range of admissible values of $s$ in the case of non-degenerate manifold established in \cite{bere,hua3,VV}, our range here is much better.  Moreover, Theorem \ref{t1} and \ref{t3}, to the best of our knowledge,  are not previously known for any single affine subspace. 

If we take $\psi(q)=q^{-\nu}$ in Theorem \ref{t3}, then we immediately obtain  an upper bound for the Hausdorff dimension of the set $\mathscr{S}(q^{-\nu})\cap\mathscr{H}$.

\begin{cor}
Let $\mathscr{H}$ be given by \eqref{e1.1} and suppose that $\omega(\bm{\widetilde{A}}) \le \frac{d+1}\nu$ for some 
$\nu\ge\frac1n$. Then we have
$$
\dim\mathscr{S}(q^{-\nu})\cap\mathscr{H}\le d-\frac{\nu n-1}{\nu+1}.
$$
\end{cor}
If the divergence counterpart of Theorem \ref{t3} is true, the above inequality actually becomes equality.

There is also a dual side of the generalized Baker-Schmidt problem. We refer the readers to \cite{BDV1, hua2, HSS} for non-degenerate manifolds and \cite{BGGV} for affine subspaces.

Throughout the paper, we will use Vinogradov's notation $f(x)\ll g(x)$ and Landau's notation $f(x)=O(g(x))$ to mean there exists a positive constant $C$ such that $|f(x)|\le Cg(x)$ as $x\to\infty$. We also denote $e^{2\pi i x}$ by $e(x)$.

\section{Multiplicatively badly approximable  matrices}\label{s2}
In this section, we state  estimates for a certain sum of products of reciprocals of fractional parts, whose proof is presented in \S \ref{s4}.

For positive integers $l, m$ and $\omega>0$, we denote the following set of multiplicatively badly approximable matrices 
$$
\left\{\bm{M} \in \R^{l \times m}\middle|\inf_{\bm{j}\in\Z^m\backslash\{\bm{0}\}}|\bm{j}|^\omega\prod_{u = 1}^l\|\bm{j}\cdot\mathrm{row}_u(\bm{M})\|>0\right\}
$$
by $\mathbf{MAD}(l, m, \omega)$.

\begin{thm}\label{t7}Let $\bm{M} \in \mathbf{MAD}(l, m, \omega)$. 
Then for all positive integers $J\ge2$ we have
$$\displaystyle \sum_{\substack{\bm{j}\in\Z^m\backslash\{\bm{0}\}\\|\bm{j}|\le J}}\prod_{u = 1}^l\|\bm{j} \cdot \mathrm{row}_u(\bm{M})\|^{-1}
\ll J^\omega (\log J)^l,$$
where the implied constant only depends on $\bm{M}$.
\end{thm}

The special case $\bm{M}\in\mathbf{MAD}(l, m, m)$ of Theorem \ref{t7} has been obtained by L\^{e} and Vaaler \cite[Theorem 2.1]{LV}. It is worth noting that their proof is based on an analytic large sieve for matrices, while ours is completely elementary. In fact, we generalize the well known gap principle (cf. \cite[Lemma 3.3]{KN}) to the matrix setting. We will prove a slightly more general version of Theorem \ref{t7}, where the notion of multiplicatively bad approximability is replaced with $\phi$-bad approximability and the $\ll$-constant is made explicit.

For a non-increasing function $\phi: \N\to(0,1)$, we say a matrix $\bm{M} \in \R^{l \times m}$ is $\phi$-badly approximable if 
$$\displaystyle\prod_{u = 1}^l\|\bm{j}\cdot\mathrm{row}_u(\bm{M})\| \ge \phi(|\bm{j}|)$$ holds for all vectors $\bm{j} \in \Z^m \backslash \{\bm{0}\}$.

\begin{thm}\label{t6} Suppose $\bm{M} \in \R^{l \times m}$ is a $\phi$-badly approximable matrix.
Then for all positive integers $J$ we have 
$$\displaystyle \sum_{\substack{\bm{j}\in\Z^m\backslash\{\bm{0}\}\\|\bm{j}|\le J}}\prod_{u = 1}^l\|\bm{j} \cdot \mathrm{row}_u(\bm{M})\|^{-1}
\le \frac{4^l}{\phi(2J)}{l+\left\lfloor \log_2\frac{1}{\phi(J)}\right\rfloor\choose l}.$$
\end{thm}
Here ${* \choose *}$ is the binomial coefficient. Note that the upper bound in Theorem \ref{t6} is of the order
$$
\frac{(\log_2\frac1{\phi(J)})^l}{\phi(2J)}.
$$

The special cases $(l,m)=(2,1)$ and $l=1$ have also been studied in \cite{F} and \cite{W}, using dynamical tools. Moreover, various
inhomogeneous variants of these sums have been treated in \cite{BHV, BL, Ch, CT}, in connection with the Littlewood conjecture. Recently, Fregoli has improved Theorem \ref{t6} a bit further by essentially removing one $\log_2\frac1{\phi(J)}$ from the upper bound \cite{F2}.

\section{Counting rational points near affine subspaces}\label{s3}
Obtaining sharp bounds for the number of rational points near a manifold lies at the core of establishing Khintchine and Jarn\'{i}k type theorems for the manifold. For the most recent advances on non-degenerate manifolds, see \cite{BVVZ1, hua1, hua3, hua5}. Our main object here is to obtain essentially optimal counting results for affine subspaces. 

Let $$
\bm{\widetilde{a}}=(q,{\bm{a}}),
$$
$$\mathcal{A}(q, \delta) := \#\left\{\bm{a}\in\Z^d \bigm| \bm{a}/q\in (0,1]^d , \left\|\bm{\widetilde{a}}\cdotp \widetilde{\bm{A}}\right\| < \delta \right\}$$
and
$$\mathcal{N}(Q, \delta) := \#\left\{\bm{\widetilde{a}}\in\N^{d + 1} \bigm| \left\|\bm{\widetilde{a}}\cdotp\widetilde{\bm{A}}\right\| < \delta, \left|\bm{\widetilde{a}}\right| \le Q\right\}.$$
Intuitively, $\mathcal{A}(q, \delta)$ counts the number of rational points $(\frac{\bm{a}}{q}, \frac{\bm{b}}{q})$ with $\frac{\bm{a}}{q}\in(0,1]^d$ that lie within
$O(\frac{\delta}{q})$ of $\mathscr{H}$. In our Theorems \ref{t4} and \ref{t5} below, we establish asymptotic formulae for $\mathcal{A}(q, \delta)$ and $\mathcal{N}(Q, \delta) $. 

\begin{thm}\label{t4} Let $\bm{A}\in\mathbf{MAD}(d, n-d,\omega)$, $q\in\N$ and $\delta\in(0,1/2]$. For any $\varepsilon>0$ we have
$$\left|\mathcal{A}(q, \delta)-(2\delta)^{n - d}q^d\right| \le \varepsilon\delta^{n - d}q^d +C_\varepsilon \delta^{n - d - \omega}\left(\log{\frac1\delta}\right)^d,$$
where $C_\varepsilon$ is only dependent on $\varepsilon$ and $\bm{A}$.
\end{thm}

\begin{thm}\label{t5} Let $\bm{\widetilde{A}}\in\mathbf{MAD}(d+1, n-d,\omega)$, $q\in\N$ and $\delta\in(0,1/2]$. For any $\varepsilon>0$ we have
$$\left|\mathcal{N}(Q, \delta)-(2\delta)^{n - d}{Q^{d + 1}}\right| \le \varepsilon\delta^{n - d}Q^{d+1} +C_\varepsilon \delta^{n - d - \omega}\left(\log{\frac1\delta}\right)^{d+1},$$
where $C_\varepsilon$ is only dependent on $\varepsilon$ and $\bm{\widetilde{A}}$.
\end{thm}

We believe that Theorems \ref{t4} and \ref{t5} are optimal except for the powers of $\log \frac1\delta$, but we do not attempt to justify this in the current paper. The proofs of these two theorems are presented in \S\ref{s5}.

\section{Proof of Theorem \ref{t6}}\label{s4}

Let $\bm{\alpha}_u=\mathrm{row}_u(\bm{M})$.
To prove the theorem, 
it is somewhat easier to first of all deal with the fractional part function $\{\cdot\}$ in lieu of the distance to the nearest integer function $\|\cdot\|$.  

Let $P_{\bm{j}} := \big(\{\bm{j} \cdot \bm{\alpha}_1\}, \{\bm{j} \cdot \bm{\alpha}_2\}, \ldots, \{\bm{j} \cdot \bm{\alpha}_l\}\big)$ be a point in $\R^l,$ 
$$
F(\bm{y}):=\prod_{u=1}^ly_u^{-1}, \textrm{ for } \bm{y}=(y_1,\ldots,y_l)\in\R^l
$$
and
$$
\mathcal{J}:=(\Z\cap[-J,J])^m\backslash\{\bm{0}\}.
$$

We observe that any two distinct points $P_{\bm{j}_1}, P_{\bm{j}_2}$ with $\bm{j}_1,\bm{j}_2\in\J$, are separated in the sense that
\begin{equation}
\prod_{u=1}^l\big|\{\bm{j}_1 \cdot \bm{\alpha}_u\} - \{\bm{j}_2 \cdot \bm{\alpha}_u\}\big|
\ge\prod_{u=1}^l\big\|(\bm{j}_1-\bm{j}_2)\cdot\bm{\alpha}_u\big\|
\ge\phi(2J).\label{e2.1}
\end{equation}

With the goal of estimating the sum
$
\sum_{\bm{j}\in\mathcal{J}}F(P_{\bm{j}}),
$
we shall count the number of points $P_{\bm{j}}$ with $\bm{j}\in\J$ in the region
$$R(\bm{k}) = [2^{-k_1}, 2^{-k_1+1}) \times [2^{-k_2}, 2^{-k_2+1 }) \times \cdots \times [2^{-k_l}, 2^{-k_l+1 }), $$
where $\bm{k}=(k_1, k_2, \ldots, k_l) \in \N^l$. Before doing this, we find some bounds on $\bm{k}$ in order for $R(\bm{k})$ to contain some $P_{\bm{j}}$.
Since for all
$\bm{j}\in\J$
  $$\prod_{u = 1}^l\{\bm{j} \cdot \bm{\alpha}_u\}\ge \prod_{u = 1}^l\|\bm{j} \cdot \bm{\alpha}_u\| \ge\phi(|\bm{j}|)\ge \phi(J),$$ 
 we have 
\begin{equation}\label{e2.2}
2^{l-\sum_{u = 1}^l k_u} \ge \phi(J)\ge\phi(2J)
\end{equation}
when $$\{\bm{j}\in\J: P_{\bm{j}}\in R(\bm{k})\}\neq\emptyset.$$

Let 
$$
z_u=2^{-k_u}\left(\phi(2J)\prod_{u=1}^l2^{k_u}\right)^{\frac1l}, \quad 1\le u\le l.
$$
Then we have
\begin{equation}\label{e2.3}
\prod_{u=1}^lz_u=\phi(2J)
\end{equation}
and

\begin{equation}\label{e2.4}
\frac{z_u}2\le 2^{-k_u}
\end{equation}
in view of \eqref{e2.2}.

To estimate the number of points $P_{\bm{j}}$ in $R(\bm{k})$, we place hyperrectangles parallel to the coordinate axes of size $z_1\times z_2\times \ldots\times z_l$ and centered at each $P_{\bm{j}}$.

First of all, we claim that these hyperrectangles do not overlap. This is because for any distinct $P_{\bm{j}_1}, P_{\bm{j}_2}\in R(\bm{k})$ with $\bm{j}_1,\bm{j}_2\in\J$, it follows, from \eqref{e2.1} and \eqref{e2.3}, that
$$
\prod_{u=1}^l\big|\{\bm{j}_1 \cdot \bm{\alpha}_u\} - \{\bm{j}_2 \cdot \bm{\alpha}_u\}\big|\ge \prod_{u=1}^lz_u,
$$
which implies that 
$$
\big|\{\bm{j}_1 \cdot \bm{\alpha}_u\} - \{\bm{j}_2 \cdot \bm{\alpha}_u\}\big|\ge z_u
$$
holds for 
at least some $u\in\{1, 2, \ldots, l\}$. This, of course, guarantees that the hyperrectangles centered at $P_{\bm{j}_1}$ and $P_{\bm{j}_2}\in\J$ do not overlap. 

Next, we observe, from \eqref{e2.4}, that each hyperrectangle centered at $P_{\bm{j}}\in R(\bm{k})$ contributes a volume of at least $\prod_{u=1}^l\frac{z_u}2$ to $R(\bm{k})$ and the total volume of $R(\bm{k})$ is $\prod_{u=1}^l2^{-k_u}$. Hence
we deduce that
$$
\#\{\bm{j}\in\J: P_{\bm{j}}\in R(\bm{k})\}\le \frac{\prod_{u=1}^l2^{-k_u}}{\prod_{u=1}^l\frac{z_u}2}
\overset{\eqref{e2.3}}{=} \frac{2^l\prod_{u=1}^l2^{-k_u}}{\phi(2J)}.
$$
Moreover, for each $P_{\bm{j}}\in R(\bm{k})$
$$
F(P_j)\le \prod_{u=1}^l2^{k_u}.
$$
It then follows that
\begin{align*}
\sum_{\substack{\bm{j}\in\J\\P_{\bm{j}}\in R(\bm{k})}}F(P_{\bm{j}})&\le 
\#\{\bm{j}\in\J: P_{\bm{j}}\in R(\bm{k})\}\prod_{u=1}^l2^{k_u}\\
&\le \frac{2^l}{\phi(2J)}.
\end{align*}

On the other hand, it follows from \eqref{e2.2} that $$\sum_{u = 1}^l k_u \le l+ \log_2\frac1{\phi(J)},$$
for which there are 
$$
{l+ \lfloor\log_2\frac1{\phi(J)}\rfloor \choose l}
$$
solutions in $\bm{k}\in\N^l$.

Hence, combining the above estimates yields
\begin{align*}
\sum_{\bm{j}\in\J}F(P_{\bm{j}})&=\sum_{\bm{k}\in\N^l}\sum_{\substack{\bm{j}\in\J\\P_{\bm{j}}\in R(\bm{k})}}F(P_{\bm{j}})\\
&\le \frac{2^l}{\phi(2J)}{l+ \lfloor\log_2\frac1{\phi(J)}\rfloor \choose l}.
\end{align*}

Now we observe that in view of $\|y\|=\min(\{y\},\{-y\})$
$$
\frac1{\|y\|}<\frac1{\{y\}}+\frac1{\{-y\}}.
$$
Therefore
\begin{align*}
&\sum_{\bm{j}\in\J}\prod_{u = 1}^l\|\bm{j} \cdot \bm{\alpha_u}\|^{-1}\\
<&\sum_{\bm{j}\in\J}\prod_{u = 1}^l\left(\frac1{\{\bm{j}\cdot\bm{\alpha}_u\}}+\frac1{\{\bm{j}\cdot(-\bm{\alpha}_u)\}}\right)\\
\le&\sum_{\bm{b}\in\{-1,1\}^l}\sum_{\bm{j}\in\J}F\big(\{b_1\bm{\alpha}_1\cdot\bm{j}\},\ldots, \{b_l\bm{\alpha}_l\cdot\bm{j}\}\big)\\
\le&\frac{4^l}{\phi(2J)}{l+ \lfloor\log_2\frac1{\phi(J)}\rfloor \choose l},
\end{align*}
where in the last inequality we use the fact that when $(\bm{\alpha}_1,\ldots,\bm{\alpha}_l)$ is $\phi$-badly approximable, so are $(b_1\bm{\alpha}_1,\ldots,b_l\bm{\alpha}_l)$ for $\bm{b}\in\{-1,1\}^l$.

\section{Proof of Theorems \ref{t4} and \ref{t5}}\label{s5}
We recall some basic properties of Selberg's magic functions. See \cite[Chapter 1]{mo} for details about the construction of these functions.

Let $\Delta=(-\delta, \delta)$ be an arc of $\R/\Z$, and $\chi_{_\Delta}(y)$ be its  characteristic function. Then there exist finite trigonometric polynomials of degree at most $J$
$$S^{\pm}_J(y)=\sum_{|j|\le J}b_j^{\pm}e(jy)$$ such that 
$$
S^{-}_J(y)\le \chi_{_\Delta}(y)\le S^{+}_J(y),
$$
$$b_0^\pm=2\delta\pm\frac1{J+1}
$$
and
$$|b_j^\pm|\le\frac1{J+1}+\min\left(2\delta,\frac1{\pi|j|}\right)
$$
for $0<|j|\le J$.  

  Therefore we have
\begin{align*}
\mathcal{A}(q, \delta) &\le \sum_{\bm{a}\in(\Z\cap[1,q])^d}\prod_{v = 1}^{n-d}S_J^+({\widetilde{\bm{a}}}\cdotp\mathrm{col}_v(\widetilde{\bm{A}}))\\
&= \sum_{\bm{a}\in(\Z\cap[1,q])^d}\prod_{v = 1}^{n-d}\sum_{|j_v| \le J}b_{j_v}^+
e(j_v{\widetilde{\bm{a}}}\cdot\mathrm{col}_v(\bm{\widetilde{A}}))\\
&= \sum_{\bm{a}\in(\Z\cap[1,q])^d}\sum_{\substack{\bm{j}\in\Z^{n-d}\\|\bm{j}| \le J}}\left(\prod_{v = 1}^{n-d}b_{j_v}^+\right)
e\left(\sum_{v = 1}^{n-d}j_v{\widetilde{\bm{a}}}\cdotp\mathrm{col}_v(\bm{\widetilde{A}})\right)\\
&= {q^d}{b_0^+}^{n-d} + \sum_{0 < |\bm{j}| \le J} \left(\prod_{v = 1}^{n-d}b_{j_v}^+\right)
\sum_{\bm{a}\in(\Z\cap[1,q])^d}e\left(\widetilde{\bm{a}}\cdotp\bm{\widetilde{A}}\cdotp\bm{j}^\intercal\right)\\
&\le \left(2\delta+\frac1{J+1}\right)^{n-d}\left(q^d+\sum_{0 < |\bm{j}| \le
     J}\prod_{u = 1}^d\|\bm{j}\cdot\mathrm{row}_u(\bm{A})\|^{-1}\right)\\
&\le \left(2\delta+\frac1{J+1}\right)^{n-d}\left(q^d+O\left(J^\omega(\log J)^d\right)\right)
\end{align*}
where the second to last line follows from the inequality $$|e(x) + e(2x) + \dots + e(qx)| \le \frac{2}{|e(x) - 1|} \le \frac{1}{\|x\|}$$ and the
last line follows by Theorem \ref{t7}.

Similarly, using the lower bound function $S_J^-(y)$ we obtain
$$
\mathcal{A}(q,\delta)\ge\left(2\delta-\frac1{J+1}\right)^{n-d}q^d+O\left(\left(2\delta+\frac1{J+1}\right)^{n-d}J^\omega(\log J)^d\right).
$$

As $J$ is a free parameter, we can set $J=\frac\kappa\delta$ for some large $\kappa>0$. Therefore, 
\begin{align*}
&\left|\mathcal{A}(q, \delta)-(2\delta)^{n - d}q^d\right| \\
\le& \left|(2\pm1/\kappa)^{n-d}-2^{n-d}\right|\delta^{n - d}q^d +O_\kappa\left( \delta^{n - d - \omega}\left(\log{\frac\kappa\delta}\right)^d\right).
\end{align*}
Now, after choosing $\kappa$ to be sufficiently large,
we complete the proof of Theorem \ref{t4}.

We may treat $\mathcal{N}(Q,\delta)$  in exactly the same manner. One only needs to note that
$$
\left|\sum_{\widetilde{\bm{a}}\in(\Z\cap[1,Q])^{d+1}}e\left(\widetilde{\bm{a}}\cdotp\bm{\widetilde{A}}\cdotp\bm{j}^\intercal\right)\right|\le \prod_{u = 1}^{d+1}\left\|\bm{j}\cdot\mathrm{row}_u(\bm{\widetilde{A}})\right\|^{-1}
$$
and the rest of the argument works without modification.

\section{Proof of Theorems \ref{t1} and \ref{t3}}\label{s6}
We state a lemma first.
\begin{lem}[Hausdorff-Cantelli, \cite{BD}]\label{l1}
Let $\{H_k\}$ be a sequence of hypercubes in $\R^n$ with side lengths $\ell(H_k)$ and suppose that for some $s>0$,
$$
\sum_k\ell(H_k)^s<\infty,
$$
then $\mathcal{H}^s(\limsup H_k)=0$.
\end{lem}
Theorems \ref{t1} and \ref{t3} are trivially true when $s>d$, so we may assume $s\le d$. Also, due to the measure theoretic nature,
 it suffices to prove Theorems $\ref{t1}$ and $\ref{t3}$ by restricting  $\bm{x}\in(0, 1]^d$.

Let $I=(0,1]^d$ and 
 $\Omega(\bm{\widetilde{A}}, \psi)$ denote the projection of $\mathscr{H}\cap \mathscr{S}_n(\psi)$ onto $I$. Since the projection map $(\bm{x},\widetilde{\bm{x}}\cdot\bm{\widetilde{A}})\to\bm{x}$ is bi-Lipschitz, we know that
$$\mathcal{H}^s(\mathscr{H}\cap \mathscr{S}_n(\psi))=0\quad\Longleftrightarrow \quad\mathcal{H}^s(\Omega(\bm{\widetilde{A}}, \psi))=0.$$ 

Let $\sigma(\bm{p},q)$ denote the set of $\bm{x}\in I$ satisfying the system of inequalities
$$
\left\{
\begin{array}{l}
|\bm{x}-\bm{a}/q|<\psi(q)/q\\
|\bm{\widetilde{x}}\cdot\bm{\widetilde{A}}-\bm{b}/q|<\psi(q)/q
\end{array}
\right.
$$ where $\bm{p}=(\bm{a},\bm{b})\in\N^d\times\Z^{n-d}$. Clearly 
$$
\Omega(\bm{\widetilde{A}},\psi)=\limsup_{(\bm{p},q)\in\N^d\times\Z^{n-d}\times\N}\sigma(\bm{p},q).
$$
Moreover when $\sigma(\bm{p},q)\neq \emptyset$, we have
$$
\left|(1,\bm{a}/q)\cdotp\bm{\widetilde{A}}-\frac{\bm{b}}q\right|\le\left|(1,\bm{a}/q)\cdotp\bm{\widetilde{A}}-\bm{\widetilde{x}}\cdot\bm{\widetilde{A}}\right|+\left|\bm{\widetilde{x}}\cdot\bm{\widetilde{A}}-\frac{\bm{b}}q\right|< C\frac{\psi(q)}{q}
$$
where $C>1$ depends only on $\bm{A}$. Then it follows that
$$
\left\|\bm{\widetilde{a}}\cdotp \widetilde{\bm{A}}\right\| < C\psi(q).
$$

\subsection{Proof of Theorem \ref{t1}}\label{s4.1}
Now to prove Theorem \ref{t1}, it suffices to show that  $\mathcal{H}^s(\Omega(\bm{\widetilde{A}},\psi))=0$ under the conditions 
\begin{equation}\label{e4.1}
\sum_{q = 1}^\infty\psi(q)^{n - d + s}q^{d - s}<\infty
\end{equation}
and
\begin{equation}\label{e4.2}
\omega(\bm{A}) < \frac{d(n - d + s)}{d + 1 - s}.
\end{equation}

In view of \eqref{e4.2}, we can choose a sufficiently small $\varepsilon>0$ such that
$\eta = \frac{d}{ \omega(\bm{A})+ \varepsilon} > \frac{d + 1 - s}{n - d + s}$. It is easily verified that there is no loss of generality in assuming that 
\begin{equation}\label{e4.3}
\psi(q)\ge q^{-\eta} \quad\text{for all }q\ge1,
\end{equation}
because if \eqref{e4.3} fails, one may replace $\psi$ with $\hat\psi(q):=\max(\psi(q), q^{-\eta})$ which satisfies \eqref{e4.1}, and then use the fact that $\mathscr{S}_n(\psi)\subseteq\mathscr{S}_n(\hat\psi)$.

Also by \eqref{e4.1}, $\psi(q)\to0$ as $q\to\infty$, hence we may assume that $\psi(q)\le \frac1{2C}$ for all $q\ge q_0$. Therefore, for each $q\ge q_0$, there are at most $\mathcal{A}(q,C\psi(q))$ nonempty $\sigma(\bm{p},q)$, each of which is contained in a hypercube with side length $\frac{2\psi(q)}q$.

Finally, 
on noting that $\bm{A}\in\mathbf{MAD}(d,n-d,\omega(\bm{A})+\frac\varepsilon2)$ by definition
and $(\log\frac1\delta)^d\ll\delta^{-\frac\varepsilon2}$, we have,
by Theorem \ref{t4}, that
$$
\mathcal{A}(q,C\psi(q))\ll \psi(q)^{n-d}q^d+\psi(q)^{n-d-\omega(\bm{A})-\varepsilon}\overset{\eqref{e4.3}}{\ll} \psi(q)^{n-d}q^d.
$$
Therefore
\begin{align*}
&\sum_{q=1}^\infty \mathcal{A}(q,C\psi(q))\left(\frac{2\psi(q)}q\right)^s\\
\ll&\sum_{q = 1}^\infty\psi(q)^{n-d}q^d\left(\frac{\psi(q)}q\right)^s\\
\ll&\sum_{q = 1}^\infty\psi(q)^{n - d + s}q^{d - s}\\
\overset{\eqref{e4.1}}{<}&\infty.
\end{align*}
Then Theorem \ref{t1} follows by  the Hausdorff-Cantelli lemma.

\subsection{Proof of Theorem \ref{t3}}
Since $\psi$ is non-increasing, it follows from \eqref{e4.1} that
\begin{equation}\label{e4.4}
\sum_{k = 1}^\infty\psi(2^k)^{n - d + s}2^{k(d - s+1)}<\infty.
\end{equation}
By the same argument as in \S\ref{s4.1}, we may assume for some $\varepsilon>0$ and $\eta= \frac{d+1}{ \omega(\bm{\widetilde{A}})+ \varepsilon} > \frac{d + 1 - s}{n - d + s}$ that
\begin{equation}\label{e4.5}
\psi(q)\ge q^{-\eta} \quad\text{for all }q\ge1.
\end{equation}
Also, there are at most $\mathcal{N}(2^{k+1},C\psi(2^k))$ nonempty $\sigma(\bm{p},q)$, for $q_0\le2^k\le q<2^{k+1}$, each of which is contained in a hypercube with side length $\frac{2\psi(2^k)}{2^k}$.

Now, Theorem \ref{t5}  implies that
\begin{align*}
\mathcal{N}(2^{k+1},C\psi(2^k))\ll &\psi(2^k)^{n-d}2^{k(d+1)}+\psi(2^k)^{n-d-\omega(\bm{\widetilde{A}})-\varepsilon}\\
\overset{\eqref{e4.5}}{\ll}& \psi(2^k)^{n-d}2^{k(d+1)}.
\end{align*}

Therefore 
\begin{align*}
&\sum_{k=1}^\infty \mathcal{N}(2^{k+1},C\psi(2^k))\left(\frac{2\psi(2^k)}{2^k}\right)^s\\
\ll&\sum_{k = 1}^\infty\psi(2^k)^{n-d}2^{k(d+1)}\left(\frac{\psi(2^k)}{2^k}\right)^s\\
\ll&\sum_{k = 1}^\infty\psi(2^k)^{n - d + s}2^{k(d - s+1)}\\
\overset{\eqref{e4.4}}{<}&\infty.
\end{align*}
Then Theorem \ref{t3} follows by  the Hausdorff-Cantelli lemma.

\proof[Acknowledgments]
We are grateful to the anonymous referees for many helpful suggestions and comments.  We also
 would like to thank Felipe Ram\'{i}rez for pointing out the reference \cite[Theorem 5]{R}.

\end{document}